\documentclass[final,leqno,onefignum,onetabnum]{siamltex1213}

\usepackage{graphicx, epstopdf}

\usepackage{amssymb}
\usepackage{amsmath, amsfonts}
\newcommand {\RR} {\mathbb R}

\newtheorem {rem} {Remark}
\newtheorem {prop} {Proposition}
\newtheorem {cor} {Corollary}

\newtheorem {defin} {Definition}

\newcommand {\ttbox} [1] {\mbox {\ttfamily*1}} 

\title {
Optimal recovery of integral operators and its applications\thanks{This work was supported by Simons Collaboration Grant N. 210363}}

\author{V. Babenko\footnotemark[2]\ 
\and Y. Babenko\footnotemark[3]\ 
\and N. Parfinovych\footnotemark[2]\
\and D. Skorokhodov\footnotemark[2]\
}

\begin {document}

 \maketitle
\slugger{sima}{xxxx}{xx}{x}{x--x}

\footnotetext[2]{Department of Mathematics and Mechanics, Dnepropetrovsk National University, Gagarina pr., 72, Dnepropetrovsk, 49010, UKRAINE}
\footnotetext[3]{Department of Mathematics, Kennesaw State University, 1100 South Marietta Pkwy, MD \# 9085, Marietta, GA, 30060, USA}

\begin{abstract}
In this paper we present the solution to the problem of recovering rather arbitrary integral operator based on incomplete information with error. We apply the main result to obtain optimal methods of recovery and compute the optimal error for the solutions to certain integral equations as well as boundary and initial value problems for various PDE's.
\end{abstract}



\pagestyle{myheadings}
\thispagestyle{plain}
\markboth{Babenko V., Babenko Y., Parfinovych, Skorokhodov}{Optimal recovery of integral operators and its applications}

\section{Introduction} \label{Intro}

Solutions to boundary (or initial) value problems for various partial differential equations require knowledge of a boundary (or initial) function. However, often time, those functions are not fully known and only partial information about them can be measured, e.g. values at some finite set of points, average values over small measurement intervals, values of $N$ first consecutive Fourier coefficients, etc. Thus, it is very important to find an approximate solution based on available information on the boundary (or initial) function. Furthermore, it is also natural and important to develop methods that provide an optimal (in some sense) approximation to the true solution. These research questions have been explored under the theory of optimal recovery of functions and operators, which is an area of Approximation Theory that started to develop in 1970s. More information on the development of the area can be found, for instance, in ~\cite{Micc_Rivl, Traub_Wozn_80, Micc_Riv_85, Korn_86, MIO_91, Plaskota, KYO, Korn}.

As for specific applications to recovering solutions of boundary and initial value problems, Magaril-Ill'yaev, Osipenko, and co-authors (see, for instance, ~\cite{MIO, Vysk_Osip_07, MIO_09}) have considered the problem of optimal $L_2$-approximation of the solution to the Dirichlet problem for Laplace's and Possion's equations in simple domains (disk, ball, annulus) based on the first $N$ consecutive Fourier coefficients of the boundary function (possibly given with an error). In order to solve this problem they have used methods of Harmonic Analysis and general results from Optimization Theory.

In this paper we address related questions of optimal approximation of the solution to several types of integral equations, boundary and initial value problems for PDE's.
We begin by solving a more general problem of recovering a rather arbitrary integral operator and sum of operators. We then present the optimal {\it method of recovery} as well as the {\it optimal error}. Next, we apply this general result to recover solutions to various boundary and initial value problems.
Moreover, we present optimal methods of recovery of the solution to boundary-value problems based on this incomplete information with error.
Naturally, the solution to the problem when information with error is used will also lead to the solution to the problem with exact information.
In this paper we focus on considering the Volterra's and Fredholm's linear integral equations as well as boundary value problems for wave, heat, and Poisson's equations. Nevertheless, the developed method is more general and can be applied to other similar problems.

The paper is organized as follows. Section \ref{S2} contains necessary definitions and notation as well as the formulation and solution of the main problem. In Section \ref{S3}, we solve the problem of optimal recovery of positive integral operators on classes of functions defined by moduli of continuity, based on information with an error about values of such functions at a fixed system of points. In Section \ref{S4}, we use our general result from Section \ref{S3} to address optimal recovery problems for the solutions of Volterra and Fredholm integral equations of the second kind, systems of linear first order differential equations with constant coefficients, Poisson's equation, the heat and wave equations.

\section{Statement and solution of the main problem} \label{S2}

\subsection{Definitions and notation}

For $l,m\in\mathbb{N}$, we let $\left\{X_j\right\}_{j=1}^{m}$ be a collection of real linear spaces, $\left\{Y_i\right\}_{i=1}^l$ be a collection of real linear normed spaces, and $\left\{Z_j\right\}_{j=1}^m$ be a collection of real linear spaces. Set 
\[
	\overline{X}: = X_1\times\ldots\times X_m, \qquad \overline{Y}:=Y_1\times\ldots\times Y_l,\qquad 	\overline{Z}:=Z_1\times\ldots\times Z_m.
\]
We write elements of spaces $\overline{X}$, $\overline{Y}$, $\overline{Z}$ as vector-columns, {\it e.g.} $\overline{x}\in\overline{X}$ is a vector-column consisting of elements $x_1,x_2,\ldots,x_m$ with $x_j\in X_j$, $j={1,\dots, m}$. This allows us to equip spaces $\overline{X}$, $\overline{Y}$, $\overline{Z}$ with natural coordinate-wise linear structure. In addition, in the space $\overline{Y}$ we introduce the norm
\begin{equation}
\label{psi_norm}
	\|\overline{y}\|_{\overline{Y}} = \|\overline{y}\|_\psi := \psi\left( \|y_1\|_{Y_1},\ldots ,\|y_l\|_{Y_l}\right),
\end{equation}
where $\psi$ is an arbitrary norm in $\mathbb{R}^l$, monotone with respect to the natural partial order in $\mathbb{R}^l$. 

By $\theta$ we denote zero of a linear space. It will be clear from the context what space is being discussed and, hence, we omit specifying it in the notation. 


Next, for a collection of linear operators $A_{ij}:X_j\to Y_i$, $i={1,\dots,l}$, and $j={1,\dots, m}$, with domains of definition $D\left(A_{ij}\right)$ we consider operator matrix
\[
	\overline{A} := \left(\!\!\!\begin{array}{cccc}
		A_{11} & A_{12} & \cdots & A_{1m} \\
		A_{21} & A_{22} & \cdots & A_{2m} \\
		\vdots & \vdots & \ddots & \vdots \\
		A_{l1} & A_{l2} & \cdots & A_{lm}
	\end{array}\!\!\!\right).
\]
The matrix $\overline{A}$ defines the operator $\overline{A}:\overline{X}\to\overline{Y}$ mapping an element $\overline{x}\in\overline{X}$ into the element $\overline{y} = \overline{A}\overline{x}$, which is a result of formal multiplication of matrix $\overline{A}$ by the vector-column $\overline{x}$, {\it i.e.}, for every $i=1,\ldots,l$, the element $y_i$ is defined as
\[
	y_i = \sum\limits_{j=1}^m A_{ij}x_j.
\]

For a given set of numbers $\sigma_j\in\{-1,1\}$, $j={1,\dots,m}$, by $\sigma$ we denote the diagonal matrix
\[
	\sigma = \textrm{diag}\,(\sigma_1,\ldots,\sigma_m) = \left(\!\!\!\begin{array}{ccccc}
	\sigma_1 & 0 & \cdots & 0 \\
	0 & \sigma_2 & \cdots & 0 \\
	\vdots & \vdots & \ddots & \vdots\\
	0 & 0 & \cdots  & \sigma_m
	\end{array}\!\!\!\right).
\]
We define the product of operator matrix $\overline{A}$ by matrix $\sigma$ as a result of formal multiplication of corresponding matrices, {\it i.e.} $\overline{A}\sigma$ is the operator matrix:
\[
	\overline{A}\sigma = \left(\!\!\!\begin{array}{cccc}
	\sigma_1A_{11} & \sigma_2A_{12} & \cdots & \sigma_mA_{1m} \\
	\sigma_1A_{21} & \sigma_2A_{22} & \cdots & \sigma_mA_{2m} \\
	\vdots & \vdots & \ddots & \vdots \\
	\sigma_1A_{l1} & \sigma_2A_{l2} & \cdots & \sigma_mA_{lm}
	\end{array}\right).
\]

In each of spaces $X_j$, $j={1,\dots,m}$, we select a class of elements $W_j\subset \bigcap\limits_{i=1}^l D(A_{ij})$, and consider the Cartesian products of classes $W_j$'s and their linear spans, respectively:
\[
	\overline{W} := W_1\times\ldots\times W_m\subset \overline{X},\quad \textrm{and}\quad \overline{\textrm{span}\,W} := \textrm{span}\,W_1\times\ldots \times\textrm{span}\,W_m.
\]
Let us assume that a collection of operators $I_j:\textrm{span}\,W_j\to Z_j$, $j={1,\dots,m}$, is given. We call {\it information operators} each of operators $I_j$ and the operator matrix $\overline{I} = \textrm{diag}\left(I_1,\ldots,I_m\right):\overline{\textrm{span}\,W}\to\overline{Z}$ as well. Note that in majority of applications we consider in this paper, information operators will be linear.

Finally, for a collection of operators $L_j:Z_j\to X_j$, $j={1,\dots,m}$, by $\overline{A}\,\overline{L}$ we denote  the operator matrix obtained as a result of formal multiplication of operator matrices $\overline{A}$ and $\overline{L} = \textrm{diag}\,(L_1,\ldots,L_m)$, where by ${A}_{ij}L_j$ we understand the composition of operators $A_{ij}$ and $L_j$.

\subsection{Optimal recovery problem and general lower estimate for the error of recovery}

In this paper we consider the problem of optimal recovery of operator $\overline{A}$ on the class $\overline{W}$ using information $\overline{I}\overline{x}$ on elements $\overline{x}$ from this class. 

For non-empty sets $U_1,\ldots U_m$ in spaces $Z_1,\ldots,Z_m$, respectively, we set 
\[
	\overline{U} = U_1\times\ldots\times U_m.
\] 
We assume that instead of $\overline{I}\overline{x}$, we know some element $\overline{z} \in \overline{I}\overline{x} + \overline{U}$, where $\overline U$ is given set, containing zero $\theta$. In this situation we say that information $\overline{I}\overline{x}$ {\it is known with $\overline{U}$-error}. 
Note that if $\overline{U}$ coincides with the origin of $\overline{Z}$ then $\overline{z} = \overline{I}\overline{x}$, and we say that information $\overline{I}\overline{x}$ is given exactly.

An arbitrary mapping $\Phi:\overline{Z}\to\overline{Y}$ is called a {\it method of recovery}. Given the operator $\overline{A}$, class $\overline{W}$, information $\overline{I}$ with $\overline{U}$-error, we define {\it the error of recovery} of operator $\overline{A}$ with the help of method $\Phi$ as follows
\[
	\mathcal{E}\left(\overline{A};\overline{W};\overline{I};{\overline{U}};\Phi\right) = \mathcal{E}_\psi\left(\overline{A};\overline{W};\overline{I};{\overline{U}};\Phi\right) := \sup\limits_{\overline{x}\in \overline{W}} \;\; \sup\limits_{\overline{z}\in \overline{I}\overline{x}+\overline{U}} \left\|\overline{A}\overline{x} - \Phi\overline{z}\right\|_\psi,
\]
and {\it the error of optimal recovery} of operator $\overline{A}$ as
\begin{equation}
\label{error_of_recovery_definition}
	\mathcal{E}\left(\overline{A};\overline{W};\overline{I};\overline{U}\right) = \mathcal{E}_\psi\left(\overline{A};\overline{W};\overline{I};\overline{U}\right) := \inf\limits_{\Phi:\overline{Z}\to \overline{Y}} \mathcal{E}\left(\overline{A};\overline{W};\overline{I};{\overline{U}};\Phi\right).
\end{equation}

{\bf The problem of optimal recovery of operator $\overline{A}$}: find the optimal error $\mathcal{E}\left(\overline{A};\overline{W};\overline{I};\overline{U}\right)$ and the method of recovery $\Phi$ (if any exists) delivering the $\inf$ in the right hand part of~(\ref{error_of_recovery_definition}). 

Clearly, when $\overline{U}$ is the origin of $\overline{Z}$, problem~(\ref{error_of_recovery_definition}) reduces to the problem of optimal recovery of operator $\overline{A}$ on the class $\overline{W}$ based on exact information $\overline{I}\overline{x}$ on elements $\overline{x}\in \overline{W}$.

Problems of optimal recovery of operators based on exact information were studied in~\cite{Smolyak,Bakhvalov,Golomb,Micc_Rivl}, and on approximate information in~\cite{March_Osip,Melk_Micc,Micc_Riv_85,MIO_91,MIO,BBP}. We also refer the reader to the discussion of closely related questions in~\cite{Traub_Wozn_80,Micc_Rivl,Korn_86,Zhensykbaev,Plaskota,Arest_Gab_95,Arestov_96}.

Let us provide the lower estimate on the error $\mathcal{E}\left(\overline{A};\overline{W};\overline{I};\overline{U}\right)$. By $\overline W (\overline U)$ we denote the following class:
$$
\overline W (\overline U):=\left\{\overline{x}\in\overline{W}:\;\;\; \left(\overline{I}\overline{x} + \overline{U}\right)\cap \left(\overline{I}(-\overline{x}) + \overline{U}\right) \ne \emptyset\right\}
$$

\begin{prop}
	\label{Prop:Lower_estimate_for_recovery_error_general}
	If for every $i=1,\dots, l$ and $j=1,\dots, m$, the operator $A_{ij}$ is odd, $\overline W (\overline U)\ne \emptyset$, and the class ${W}_j$ is centrally symmetric, then
	\begin{equation}
	\label{lower_estimate_general}
	\mathcal{E}_\psi\left(\overline{A};\overline{W};\overline{I};\overline{U}\right)  \geqslant \sup\limits_{\overline{x}\in\overline{W} (\overline U)} \left\|\overline{A}\overline{x}\right\|_\psi .
	\end{equation}
\end{prop}

\begin{proof}
	For any method of recovery $\Phi:\overline{Z}\to \overline{Y}$, 
	\[
		\begin{array}{l}
	\mathcal{E}_\psi\left(\overline{A};\overline{W};\overline{I};{\overline{U}};\Phi\right) =  \sup\limits_{\left\{\overline{x}\in\overline{W}\atop \overline{z}\in \overline{I}\overline{x}+{\overline{U}}\right\}} \|\overline{A}\overline{x}- \Phi\overline{z}\|_\psi\\
	\qquad\quad\geqslant \displaystyle \max{\left\{\sup\limits_{\left\{\overline{x}\in\overline{W}\atop \overline{z}\in \overline{I}\overline{x}+{\overline{U}}\right\}} \|\overline{A}\overline{x}- \Phi\overline{z}\|_\psi; \sup\limits_{\left\{-\overline{x}\in\overline{W}\atop \overline{z}\in \overline{I}(-\overline{x})+{\overline{U}}\right\}} \|-\overline{A}\overline{x}- \Phi\overline{z}\|_\psi\right\}}\\
	\qquad\quad\displaystyle \geqslant \frac 12\sup\limits_{\overline{x}\in\overline{W} (\overline U)}   \left\{\left\|\overline{A}\overline{x}-\Phi \overline{z}\right\|_\psi+ \left\|\overline{A}\overline{x}+\Phi\overline{z}\right\|_\psi\right\}\nonumber\\
	\qquad\quad \geqslant \sup\limits_{\overline{x}\in\overline{W} (\overline U)} \left\|\overline{A}\overline{x}\right\|_\psi,\nonumber
		\end{array}
	\]
	which completes the proof.
\end{proof}

In particular, when $\overline{I}$ is an even operator, the condition 
\[
	\left(\overline{I}\overline{x} + \overline{U}\right) \cap \left(\overline{I}(-\overline{x}) + \overline{U}\right)\ne \emptyset
\]
is valid for every $\overline{x}\in \overline{W}$. Also, one can easily verify that there holds the following consequence from Proposition~\ref{Prop:Lower_estimate_for_recovery_error_general}.

\begin{prop}
\label{Prop:Lower_estimate_for_recovery_error}
Let assumptions of Proposition~\ref{Prop:Lower_estimate_for_recovery_error_general} hold. In addition, for every $j=1,\dots, m$, we let $I_j$ be odd, and $U_j$ be centrally symmetric. Then
\begin{equation}
\label{lower_estimate}
	\mathcal{E}_\psi\left(\overline{A};\overline{W};\overline{I};\overline{U}\right)  \geqslant \sup\limits_{\left\{\overline{x}\in\overline{W}\atop \overline{I}\overline{x}\in  \overline{U}\right\}} \left\|\overline{A}\overline{x}\right\|_\psi .
\end{equation}
\end{prop}

Note that the lower estimate provided by Proposition~\ref{Prop:Lower_estimate_for_recovery_error} might not be sharp. Therefore, in the rest of this section we consider some general situations when inequality~(\ref{lower_estimate}) turns into equality.

\subsection{Normed lattices and positive operators. }
\label{subsec:BanachLatticesAndPositiveOperators}

Let us follow~\cite{Schaefer} in order to introduce the concepts of an ordered vector space, a normed lattice, and a positive operator.

\begin{defin}
\label{OrderedVectorSpace}
Given a linear space $X$ over the field of real numbers $\mathbb{R}$ and a partial order ``$\prec_X$'' on the set $X$, we call the pair $(X,\prec_X)$ an ordered vector space if:
\begin{enumerate}
\item $x\prec_X y$ implies $x+z \prec_X y+z$, for all $x,y,z\in X$;
\item $x\prec_X y$ implies $\lambda x \prec_X \lambda y$, for all $x,y\in X$ and $\lambda\in\mathbb{R}_+$.
\end{enumerate}
\end{defin}

In what follows, for brevity (when it does not lead to confusion), we omit mentioning partial order ``$\prec_X$'' in the notation of an ordered vector space $\left(X,\prec_X\right)$.

Also, we reserve notation ``$\leqslant$'' for the standard linear order in $\mathbb{R}$.

\begin{defin}
\label{VectorLattice}
An ordered vector space $X$ is called a vector lattice (the Riesz space), if any two elements $x,y\in X$ have supremum $x\vee y := \sup\{x; y\}$ and infimum $x\wedge y := \inf\{x; y\}$.
\end{defin}

For a vector lattice $X$, by $|x| := x\vee (-x)$ we define absolute value of $x\in X$. In addition, we call a norm on a vector lattice $X$ a {\it lattice norm}, if $|x| \prec_X |y|$ implies $\|x\| \leqslant \|y\|$ for all $x, y \in X$.

\begin{defin}
\label{NormedLattice}
A normed lattice is a real normed space $X$ endowed with an ordering ``$\prec_X$'' such that $(X,\prec_X)$ is a vector lattice and the norm on $X$ is a lattice norm.
\end{defin}

Note that given a collection of ordered vector spaces $X_j$, $j={1,\dots,m}$, their Cartesian product $\overline{X}$ is also an ordered vector space with respect to naturally defined partial order ``$\prec_{\overline{X}}$'':
\[
	(\overline{x}'\prec_{\overline{X}} \overline{x}'')\Leftrightarrow (\forall j=1,\dots, m,\; x_j'\prec_{X_j} x''_j).
\]
Similarly, for a collection of normed lattices  $Y_i$, $i={1,\dots,l}$, their Cartesian product $\overline{Y}$ is also a normed lattice with respect to the norm $\|\cdot\|_\psi$ defined by~(\ref{psi_norm}) and partial order ``$\prec_{\overline{Y}}$''.

Finally, we define the positive operator between ordered vector spaces as follows.

\begin{defin}
Let $X$ and $Y$ be ordered vector spaces. A linear operator $T:X\to Y$ is called positive if $\theta \prec_Y Tx$ whenever $\theta \prec_X x$.
\end{defin}

\subsection{General results for positive operators}

In this subsection we present some results on optimal recovery of positive operators and, in particular, identity operator. Furthermore, we show that under certain assumptions, once we know how to recover (in an optimal way) the identity operator on each of classes $W_j$, based on information $I_j$ with $U_j$-error, we can recover (in an optimal way) any operator matrix $\overline{A}$ consisting of positive linear operators (and even operator matrix $\overline{A}\sigma$) on the class $\overline{W}$, based on information $\overline{I}$ with $\overline{U}$-error.

Let $X$ be a normed lattice, and by $id_{X}$ we denote the identity operator. We start with the problem of optimal recovery of the identity operator. Let $Z$ be a real linear space, $W\subset X$ be centrally symmetric class, $I:X\to Z$ be an odd information operator, and $U\subset Z$ be non-empty centrally symmetric set.

\begin{prop}
\label{Upper_Estimate_for_identity_Operator}
If there exist an operator $L:Z\to X$ and a function $\varphi \in W$, $I\varphi \in U$, such that for any $x\in W$ and $z\in Z$ we have
\begin{equation}
\label{reference_star}
	(z\in Ix+U)\;\Rightarrow \;   (-\varphi \prec_{X} x - Lz\prec_{X} \varphi),
\end{equation}
then operator $L$ is the optimal method of recovery of $id_X$ on the class $W$, based on information $I$ with $U$-error, and
\begin{equation}
\label{equality_for_identity}
	\mathcal{E}\left(id_X;W;I;U\right) = \mathcal{E}\left(id_X;W;I;U; L\right) = \left\|\varphi\right\|_X.
\end{equation}
\end{prop}

\begin{proof} By assumption, for every $x\in W$ and $z\in Ix+U$ we have
\[
	-\varphi \prec_X x - Lz\prec_X \varphi.
\]
Since $X$ is a normed lattice, from the latter we obtain 
\[
	\left\|x - Lz\right\|_X\leqslant \|\varphi\|_X.
\]
Hence,
\[
	\mathcal{E}\left(id_X;W;I;U\right) \leqslant \mathcal{E}\left(id_X;W;I;U;L\right) = \sup\limits_{x\in W}\sup\limits_{z\in Ix+U}\left\|x - Lz\right\|_X \leqslant \|\varphi\|_X.
\]

On the other hand $I\varphi \in U$. Hence, due to  Proposition~\ref{Prop:Lower_estimate_for_recovery_error}, we obtain
\[
	\mathcal{E}\left(id_X;W;I;U\right) \geqslant \sup\limits_{\left\{x\in W \atop Ix\in U\right\}} \|x\|_X \geqslant \|\varphi\|_X.
\]
\end{proof}

Next, we present the result on optimal recovery of positive operators, which follows from Proposition~\ref{Upper_Estimate_for_identity_Operator}.

\begin{prop}
\label{Upper_Estimate_for_Positive_Operator}
Under assumptions of Proposition~\ref{Upper_Estimate_for_identity_Operator}, let $Y$ be a normed lattice and $A:X\to Y$ be a positive linear operator with domain of definition $\mathcal{D}(A)\supset W$. Then $\Phi = AL$ is the optimal method of recovery of operator $A$ on the class $W$, based on information $I$ with $U$-error, and
\[
	\mathcal{E}\left(A;W;I;U\right) = \mathcal{E}\left(A;W;I;U; \Phi\right) = \left\|A\varphi\right\|_Y.
\]
\end{prop}

\begin{proof}
Indeed, let $x\in W$ and $z\in Ix + U$ be arbitrary. According to~(\ref{reference_star}), we have $-\varphi \prec_X x - Lz \prec_X \varphi$. Hence, due to positivity of operator $A$, we obtain
\[
	-A\varphi \prec_Y Ax - ALz \prec_Y A\varphi.
\]
Taking into account that $\|\cdot\|_Y$ is a lattice norm, we deduce that
\[
	\left\|Ax - ALz\right\|_Y \leqslant \left\|A\varphi\right\|_Y
\]
and, therefore,
\[
	\mathcal{E}\left(A;W;I;U\right) \leqslant \mathcal{E}\left(A;W;I;U;AL\right) \leqslant \|A\varphi\|_Y.
\]
The opposite inequality follows from Proposition~\ref{Prop:Lower_estimate_for_recovery_error}. 
\end{proof}

\begin{rem}
	In Proposition~\ref{Upper_Estimate_for_Positive_Operator} the condition that $X$ is a normed lattice can be relaxed to the following one: $X$ is an ordered vector space.
\end{rem} 

Finally, we present the generalization of Proposition~\ref{Upper_Estimate_for_Positive_Operator} to the case of optimal recovery of operator matrices. 

In order to state the corresponding result, we let $X_1,\ldots,X_m$ be ordered vector spaces, $Z_1,\ldots,Z_m$ be real linear spaces, and $Y_1,\ldots,Y_l$ be normed lattices, $m,l\in\mathbb{N}$. In addition, let  $A_{ij}:X_j\to Y_i$, $j = 1,\dots, m$ and $i = 1,\dots,l$, be positive linear operators with domains of definition $\mathcal{D}\left(A_{ij}\right)$, and $I_j:X_j\to Z_j$, $j = 1,\dots,m$, be odd information operators. Let also $W_j\subset \bigcap\limits_{i=1}^{l} \mathcal{D}\left(A_{ij}\right)$, $j = 1,\dots,m$, be centrally symmetric classes, and $U_j\subset Z_j$, $j = 1,\dots,m$, be centrally symmetric sets. Finally, let $\psi:\mathbb{R}^l\to\mathbb{R}$ be an arbitrary norm monotone with respect to the natural partial ordering in $\mathbb{R}^l$.

\begin{theorem}
\label{Optimality}
If for every $j=1,\dots,m$, there exist an operator $L_j:Z_j\to X_j$ and a function $\varphi_j\in W_j$, $I_j\varphi_j\in U_j$, such that for any $\overline{x}\in \overline{W}$ and $\overline{z}\in \overline{Z}$, we have
\[
	\left(\overline{z}\in \overline{I}\overline{x}+\overline{U}\right)\;\Rightarrow \;   \left(\forall j=1,\dots,m,\; -\varphi_j \prec_{X_j} x_j - L_j z_j\prec_{X_j} \varphi_j\right),
\]
then for every $\sigma = \textrm{diag}\,(\sigma_1,\ldots ,\sigma_m)$, $\sigma_j\in\{-1,1\}$, $j=1,\dots, m$, the operator $\Phi = \overline{A}\sigma \overline{L}$ is the optimal method of recovery of operator $\overline{A}\sigma$ on the class $\overline{W}$ based on information $\overline{I}$ with  $\overline{U}$-error, and, furthermore,
\[
	\mathcal{E}_\psi\left(\overline{A}\sigma;\overline{W};\overline{I};\overline{U}\right) = \mathcal{E}_\psi\left(\overline{A}\sigma;\overline{W};\overline{I};\overline{U}; \Phi\right) = \left\|\overline{A}\overline{\varphi}\right\|_\psi.
\]
\end{theorem}

\begin{proof}
Let $\overline{x}\in \overline{W}$ and $\overline{z}\in \overline{I}\overline{x}+\overline{U}$. Then for any $j=1,\dots, m$, we have
\[
	-\varphi_j \prec_{X_j} x_j - L_j z_j\prec_{X_j} \varphi_j.
\]
Since operators $A_{ij}$ are linear and positive, we obtain 
\[
	-A_{ij}\varphi_j \prec_{Y_i} A_{ij}(x_j-L_jz_j) = A_{ij}x_j - A_{ij}L_j{z_j} \prec_{Y_i} A_{ij}\varphi_j.
\]
From the latter we conclude that $\forall \sigma_j \in \{-1,1\}$
\[
	-A_{ij}\varphi_j \prec_{Y_i} \sigma_jA_{ij}x_j - \sigma_jA_{ij}L_j{z_j} \prec_{Y_i} A_{ij}\varphi_j.
\]
Summing up these inequalities over $j$, we see that
\[
	-\sum\limits_{j=1}^mA_{ij}\varphi_j \prec_{Y_i} \sum\limits_{j=1}^m\sigma_jA_{ij}x_j - \sum\limits_{j=1}^m\sigma_jA_{ij}L_j{z_j} \prec_{Y_i} \sum\limits_{j=1}^mA_{ij}\varphi_j.
\]
Therefore,
\[
	\left|\sum\limits_{j=1}^m\sigma_jA_{ij}x_j - \sum\limits_{j=1}^m\sigma_jA_{ij}L_j{z_j}\right| \prec_{Y_i} \sum\limits_{j=1}^mA_{ij}\varphi_j.
\]
Since $\left\|\cdot\right\|_{Y_i}$ is the lattice norm, from the latter we derive
\[
	\left\|\sum\limits_{j=1}^m\sigma_jA_{ij}x_j - \sum\limits_{j=1}^m\sigma_jA_{ij}L_jz_j\right\|_{Y_i} \leqslant \left\|\sum\limits_{j=1}^mA_{ij}\varphi_j\right\|_{Y_i},
\]
and since $\Phi = \overline{A}\sigma\overline{L}$, 
\[
	\left\| \overline{A}\sigma\overline{x}-\overline{A}\sigma \overline{L}\overline{z}\right\|_\psi = \left\| \overline{A}\sigma\overline{x}-\Phi\overline{z}\right\|_\psi \leqslant \left\|\overline{A}\overline{\varphi}\right\|_\psi .
\]
Hence, 
\[
	\mathcal{E}_{\psi}\left(\overline{A}\sigma;\overline{W};\overline{I};\overline{U};\Phi\right) \leqslant \mathcal{E}_{\psi}\left(\overline{A}\sigma;\overline{W};\overline{I};\overline{U}\right) \leqslant\left\|\overline{A}\overline{\varphi}\right\|_{\psi}.
\]
In order to obtain the lower estimate, we observe that $\sigma\overline{\varphi}\in \overline{W}$ and $\overline{I}(\sigma\overline{\varphi})\in \overline{U}$. Due to Proposition~\ref{Prop:Lower_estimate_for_recovery_error} we obtain
\[
	\mathcal{E}_\psi\left(\overline{A}\sigma;\overline{W};\overline{I};{\overline{U}}\right) \geqslant \sup\limits_{\left\{\overline{x}\in\overline{W}\atop \overline{I}\overline{x}\in \overline{U}\right\}} \left\|\overline{A}\sigma \overline{x}\right\|_\psi\geqslant \left\|\overline{A}\sigma (\sigma\overline{\varphi}) \right\|_\psi=\left\|\overline{A}\overline{\varphi}\right\|_\psi
\]
as $\sigma\sigma = \textrm{diag}\,(1,\ldots, 1)$ is the identity matrix. 
\end{proof}

\section{Optimal recovery of integral operators}\label{S3}

In this section we introduce the concept of an integral operator and apply Theorem~\ref{Optimality} to the problem of optimal recovery of positive integral operators on classes of functions defined by moduli of continuity, based on information with an error about values of such functions at a fixed system of points.

\subsection{Integral operators on metric spaces}
\label{subsec:IntegralOperators}

We follow~\cite{Arons_Szept} (see also~\cite{Korotkov}) to introduce the notion of integral operators on metric spaces.
First, we let $(M, \mu)$ be the space with $\sigma$-finite measure, {\it i.e.} $M$ is some set and $\mu$ is a $\sigma$-finite measure on $\sigma$-algebra $\Sigma_M$ of subsets in $M$. By $\mathfrak{M}\left(M, \mu\right)$ we denote the space of all $\mu$-measurable $\mu$-{\it a.e.} finite functions defined on $M$ (identifying $\mu$-equivalent functions as usual). The space $\mathfrak{M}\left(M, \mu\right)$ is equipped with the natural partial order ``$\prec$'': for every $x',x''\in \mathfrak{M}\left(M,\mu\right)$,
\[
	\left(x'\prec x''\right)\; \Leftrightarrow\;\left(\textrm{for}\;\mu-a.\,e.\; t\in M, \;\left|x'(t)\right| \leqslant \left|x''(t)\right|\right). 
\]
Hence, we can consider the space $\mathfrak{M}\left(M,\mu\right)$ as an ordered vector space.

\begin{defin} (see~\cite{Arons_Szept},~\cite[Ch.~1, \S 2]{Korotkov}) 
Let $(M, \mu)$ and $(N,\nu$) be spaces with $\sigma$-finite positive measures,	 $E$ be a linear manifold in $\mathfrak{M}(M, \mu)$. A linear operator $T : E \to \mathfrak{M}(N, \nu)$ is called {\it integral operator} if there exists $\nu\times\mu$-measurable function $K(s,t)$ such that for every $x\in E$,
\begin{equation}\label{integral_operator}
	Tx(s) := \int_M K(s,t)\,x(t)\,d\mu(t), \qquad \textrm{for}\; \mu - a.\,e. \; \; s\in N.
\end{equation}
The above integral is understood in the Lebesgue sense. The function $K(s,t)$ is called {\it the kernel} of operator $T$.
\end{defin}

\begin{rem}
	Clearly, an integral operator is positive if its kernel is {$\nu\times \mu$-a.e.} non-negative.
\end{rem}

Next, we let $M = M_\rho$ be a metric space endowed with the metric $\rho$. By $\Sigma_{\rho}$ we denote the Borel $\sigma$-algebra of subsets of $M_{\rho}$, {\it i.e.} the minimal $\sigma$-algebra generated by open sets in $M_{\rho}$. We consider an arbitrary non-negative $\sigma$-finite measure $\mu:\Sigma_\rho \to \mathbb{R}_+$. For convenience, if the equivalence class in $\mathfrak{M}\left(M_\rho,\mu\right)$ contains a continuous function, then we identify this function with the whole equivalence class.

By $B_\mu$ and $C_\mu$ let us also denote the sets of $\mu$-essentially bounded and $\mu$-{\it a.e.} continuous functions $x:M_\rho\to\mathbb{R}$, respectively. 

For a compact set $M'\subset M_\rho$, we let $\chi_{M'}$ stand for the {\it characteristic} (or {\it indicator}) function of the set $M'$. We consider classes
\[
	\begin{array}{lll}
		B_\mu(M') & := & \left\{x \in B_\mu\,:\,\textrm{supp}\,x\subset M'\right\},\\
		C_\mu(M') & := & \left\{x = y\cdot\chi_{M'}\,:\,y\in C_\mu\right\},\\
		\tilde{C}_\mu(M') & := & \left\{x \in C_\mu\,: \,  \textrm{supp}\, x\subset M' \right\}.
	\end{array}	
\]
By definition, $\tilde{C}_\mu(M') \subset C_\mu(M') \subset B_\mu(M')$. 

\subsection{Classes $H^{\omega}$ and generalized Voronoi cells}
\label{subsec:VoronoiCells}

We recall that a function $\omega:\mathbb{R}_+\to\mathbb{R}_+$, $\mathbb{R}_+ := [0,\infty)$, is called {a modulus of continuity} (see, for example,~\cite{Korn_76}) if $\omega(0) = 0$, $\omega$ is continuous, non-decreasing, and semi-additive function. The latter means that $\omega(t'+t'')\leqslant \omega(t') + \omega(t'')$, for every $t',t''\in\mathbb{R}_+$.

We consider the problem of optimal recovery of positive integral operators on the classes defined by a modulus of continuity $\omega$:
\[
	H^{\omega}_\mu(M') := \left\{x\in C_\mu(M'):\left|x(t') - x(t'')\right|\leqslant \omega\left(\rho(t',t'')\right),\; \forall t',t''\in M'\right\},
\]
\[
	\tilde{H}^{\omega}_\mu(M') := H_\mu^\omega(M')\cap \tilde{C}_\mu(M').
\]
In addition, we assume that information on functions $x\in H^\omega_\mu\left(M'\right)$ or $x\in \tilde{H}^\omega_\mu\left(M'\right)$ is provided by information operators $I:C_{\mu}\left(M'\right) \to \mathbb{R}^n$, $n\in\mathbb{N}$, of the form
\[
	Ix = I_{Q} x := \left(x\left(q_1\right),\ldots,x\left(q_n\right)\right),\qquad x\in C_\mu(M'),
\]
where $Q = \left\{q_j\right\}_{j=1}^n$ is a fixed set of points in $M'$, and is known with $U_e$-error, $e\in\mathbb{R}_+^n$, where
\[
	U_{e} = [-e_1,e_1]\times\ldots\times[-e_n,e_n].
\] 

Next, we construct the operator $L$ that would satisfy assumptions of  Proposition~\ref{Upper_Estimate_for_identity_Operator}. To this end, we first introduce the following two functions:
\begin{equation}
\label{f_omega}
	\tau(t) = \tau_{\omega,Q,e}(t) :=
		\left\{\begin{array}{ll}
		\min\limits_{j = 1,\dots, n} \left(e_j + \omega(\rho(t,q_j))\right), & t\in M',\\
		0, & t \in M\setminus M',
	\end{array}\right.
\end{equation}
and
\begin{equation}
\label{f_0_omega}
	\tilde{\tau}(t) = \tilde{\tau}_{\omega,Q,e}(t) := 
		\min\left\{\tau_{\omega,Q,e}(t);\, \omega(\rho(t,\partial M'))\right\},\qquad t\in M,
\end{equation}
where $\rho(t,\partial M') := \inf\limits_{s\in \partial M'} \rho(t,s)$ denotes the distance between point $t\in M$ and the boundary $\partial M'$ of $M'$. One can easily verify that $\tau \in H_\mu^\omega(M')$, $\tilde{\tau} \in \tilde{H}_\mu^\omega(M')$, and $I_Q\tau = I_Q\tilde{\tau} \in U_e$.

\begin{figure}[h!]
	\centering
	\includegraphics[angle=0, width=3in]{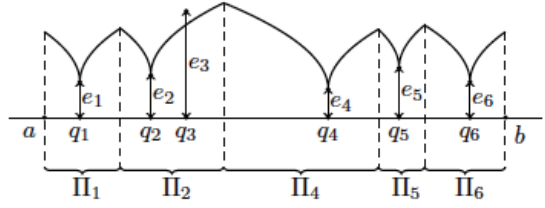}\label{pic1}
	\caption{Generalized Voronoi cells on the segment $[a,b]$}
	
\end{figure}

Next, we define {\it generalized Voronoi cells}. To this end, we first let
\[
	\tilde{\Pi}_0 = \tilde{\Pi}_0(\omega,Q,e) := \left\{t\in M'\;:\;\omega(\rho(t,\partial M')) \leqslant \tau_{\omega,Q,e}(t)\right\},
\]
and, for $j=1,\ldots,n$, we consider
\[
	\Pi'_j = \Pi'_j(\omega,Q,e):= \left\{t\in M'\,:\, \tau_{\omega,Q,e}(t) = e_{j}+ \omega(\rho\left(t,q_{j}\right)) \right\}.
\]
Then, generalized Voronoi cells (see Figure~\ref{pic1}) are defined iteratively as follows
\[
	\displaystyle\Pi_{1} := \Pi'_1,\qquad \Pi_{j} := \Pi'_{j}\setminus \bigcup\limits_{s=1}^{j-1} \Pi_s,\qquad j=2,\ldots,n,
\]
and 
\[
	\tilde{\Pi}_{j} = \tilde{\Pi}_j \setminus \tilde{\Pi}_{0},\qquad j=1,\ldots,n.
\]

Since every function $x\in C_\mu(M')$ is $\mu$-measurable on $M$ in the Borel $\sigma$-algebra $\Sigma_{\rho}$, we see that sets $\Pi_1,\Pi_2,\ldots,\Pi_n$, and sets $\tilde{\Pi}_0,\tilde{\Pi}_1,\ldots,\tilde{\Pi}_n$ are $\mu$-measurable. Moreover, due to construction, we conclude that both collections of sets are pairwise disjoint and
\[
	M' = \Pi_1\cup \ldots \cup \Pi_n = \tilde{\Pi}_0\cup\tilde{\Pi}_1\cup \ldots \cup \tilde{\Pi}_n.
\]

\begin{figure}[h!]
	\centering
	\includegraphics[angle=0, width=3in]{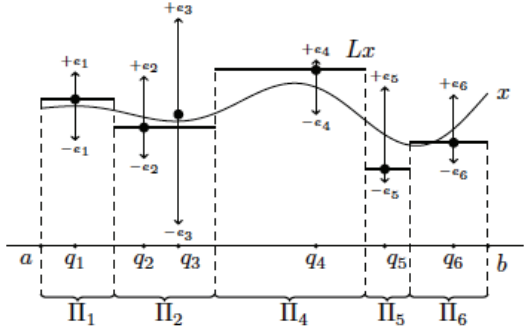}\label{pic3}
		\caption{Method of recovery of a function defined on the segment based on information about its values at 6 points with error}
\end{figure}

Finally, we define operators $L = L_{\omega,Q,e} : \mathbb{R}^n \to B_\mu\left(M'\right)$ and $\tilde{L} = \tilde{L}_{\omega,Q,e} : \mathbb{R}^n \to B_\mu\left(M'\right)$ (see Figure~\ref{pic3}) as follows: for every $z \in \mathbb{R}^n$, 
\[
	Lz(t) := \left\{\begin{array}{ll}
		z_j, & t\in \Pi_j,\;\;  j=1,\ldots,n,\\
		0, & t\in M\setminus M',
	\end{array}\right.
\]
and
\[
	\tilde{L}z(t) := \left\{\begin{array}{ll}
		z_j, & t\in \tilde{\Pi}_j, \;\; j=1,\ldots,n,\\
		0, & t\in \tilde{\Pi}_0 \cup \left(M\setminus M'\right).
	\end{array}\right.
\]

Note that when some of $e_j$'s are large enough, it is possible that some of the sets $\Pi_j$ and $\tilde{\Pi}_j$ are empty. This, in turn, means that information at the corresponding point $q_j$ is ``ignored'' by operators $L$ and $\tilde{L}$.

\subsection{Optimal recovery of operators on the class $H^\omega$}

In this section, we present some important direct corollaries from Theorem~\ref{Optimality} for positive integral operators and their sums.

For a space $\left(M,\mu\right)$ with $\sigma$-finite measure, we let 
$L_1\left(M,\mu\right)$ be the space of absolutely integrable functions $x:M\to\mathbb{R}$ with the standard norm
\[
	\|x\|_1 = \int_{M} |x(t)|\,d\mu(t).
\]

We start with the corollary from Proposition~\ref{Upper_Estimate_for_Positive_Operator} for positive integral operators. In order to state the corresponding result, we let $n\in\mathbb{N}$, $\omega$ be a modulus of continuity, $e\in\mathbb{R}_+^n$, $M_\rho$ be a metric space, $\mu$ be a $\sigma$-finite measure on $\Sigma_\rho$, $M'$ be a compact subset of $M_\rho$, $(N,\nu)$ be a space with $\sigma$-finite measure, $A:B_\mu(M')\to \mathfrak{M}\left(N,\nu\right)$ be an integral operator with the non-negative kernel $K$, $Y\subset \mathfrak{M}\left(N,\nu\right)$ be a normed lattice such that $A\left(B_\mu\left(M'\right)\right) \subset Y$, $Q = \left\{q_j\right\}_{j=1}^n\subset M'$ be a fixed system of points.

\begin{theorem}
\label{T1}
Let either $W = H^\omega_\mu(M')$, $\varphi = \tau_{\omega,Q,e}$, $L = L_{\omega,Q,e}$, or $W = \tilde{H}_\mu^\omega(M')$, $\varphi = \tilde{\tau}_{\omega,Q,e}$, $L = \tilde{L}_{\omega,Q,e}$. Then the method $\Phi = A\,L$ is the optimal method of recovery of operator $A$ on the class $W$ based on information $I_Q$ with $U_e$-error, and
\begin{equation}
\label{general_relation}
	\mathcal{E}\left(A;W;I_Q;U_e\right) = \mathcal{E}\left(A;W;I_Q;U_e;\Phi\right) = \displaystyle\left\|\, \int_{M'} K(\cdot, x)\,\varphi(x)\, d\mu(x)\right\|_{Y}.
\end{equation}
\end{theorem}

\begin{proof} We consider only the first case when $W = H_\mu^\omega(M')$, $\varphi = \tau_{\omega,Q,e}$, $L = L_{\omega,Q,e}$ as the proof of the second case follows similar arguments.

We let $x\in H_\mu^{\omega}(M')$ and $z\in Ix + U$ be given, and let $\Pi_{1},\ldots,\Pi_n$ be generalized Voronoi cells on the set $M'$. We observe that for every $j=1,\dots, n$ and $t\in \Pi_{j}$, 
\[
	\begin{array}{rcl}
		x(t) - Lz(t) & = & \displaystyle \left|x(t) - z_j\right| \leqslant \displaystyle \left|x(t)-x(q_j)\right| + \left|x(q_j)-z_j\right| \\
		& \leqslant & \displaystyle \omega(\rho(t, q_j)) + e_j = \tau_{\omega,Q,e}(t) = \varphi(t).
	\end{array}
\]
Hence, we conclude that $-\varphi \prec x - Lz \prec \varphi$ where ``$\prec$'' is the natural partial order in $\mathfrak{M}\left(M_\rho,\mu\right)$. Since $A$ is a linear positive operator, conditions of Theorem~\ref{Optimality} are satisfied, and by applying it, we complete the proof.
\end{proof}

When $Y = L_1\left(N,\nu\right)$, we obtain the following consequence from Theorem~\ref{T1}.

\begin{cor}
\label{Cor1}
Under assumptions of Theorem~\ref{T1}, we take kernel $K$ to be $\nu\times\mu$-integrable and $Y = L_1(N, \nu)$. Then
\begin{equation}
\label{L1_equality}
	\mathcal{E}\left(A;W;I_Q;U_e\right)= \int_{M'} \varphi(x) \int_{N} K(y,x)  \,d\nu(y)\, d\mu(x).
\end{equation}
In particular, when the inner integral in~(\ref{L1_equality}) is independent of $x$ and is denoted by $C_K$, we have
\[
	\mathcal{E}\left(A;W;I_Q;U_e\right)= C_K\int_{M'} \varphi(x)\,  d\mu(x).
\]
\end{cor}

Next, we state the consequence of Theorem~\ref{Optimality} for the problem of optimal recovery of sums of positive operators on classes defined by moduli of continuity. We need the following notation. Let $m,l\in\mathbb{N}$ and $\psi:\mathbb{R}^l\to\mathbb{R}$ be a norm monotone with respect to the natural partial order in $\mathbb{R}^l$. For every $j=1,\dots, m$, we let $n_j\in\mathbb{N}$; $\omega_j$ be a modulus of continuity; $e_j\in\mathbb{R}^{n_j}_+$ and $U_j = U_{e_j}$; $M_{j} = M_{\rho_j}$ be a metric space, $\mu_j$ be a $\sigma$-finite measure on $\Sigma_{\rho_j}$; $M'_j$ be a compact subset in $M_{j}$; $Q_j\subset M_j'$ be a fixed systems of $n_j$ points, and $I_j := I_{Q_j}$. For every $i=1,\dots, l$, we let $\left(N_i,\nu_i\right)$ be a space with $\sigma$-finite measure. For every $i=1,\dots, l$ and $j=1,\dots, m$, we let $A_{ij}:B_{\mu_j}\left(M_j'\right)\to \mathfrak{M}\left(N_i,\nu_i\right)$ be an integral operator with non-negative kernel $K_{ij}$, and $Y_i\subset \mathfrak{M}\left(N_i,\nu_i\right)$ be a normed lattice such that $\bigcup\limits_{j=1}^{m}A_{ij}\left(B_{\mu_j}\left(M'_j\right)\right) \subset Y_i$.

\begin{theorem}
\label{T2}
Let, for every $j=1,\dots, m$, either $W_j = H_{\mu_j}^{\omega_j}\left(M'_j\right)$, $\varphi_j = \tau_{\omega_j,Q_j,e_j}$, $L_j = L_{\omega_j,Q_j,e_j}$, or $W_j = \tilde{H}_{\mu_j}^{\omega_j}\left(M'_j\right)$, $\varphi_j = \tilde{\tau}_{\omega_j,Q_j,e_j}$, $L_{j} = \tilde{L}_{\omega_j,Q_j,e_j}$. Then, for every $\sigma =\textrm{diag}\,\left(\sigma_1,\ldots,\sigma_m\right)$, where $\sigma_j \in \{-1,1\}$, $j = 1,\dots, m$, the method $\Phi = \overline{A}\sigma\overline{L}$, where $\overline{L} = \textrm{diag}\,\left(L_1,\ldots,L_m\right)$, is the optimal method of recovery of operator $\overline{A}\sigma$ on the class $\overline{W}$ based on information $\overline{I}$ with $\overline{U}$-error, and 
\begin{equation}
\label{general_relation_sum}
	\mathcal{E}_\psi\left(\overline{A}\sigma;\overline{W};\overline{I}; \overline{U}\right) =  \mathcal{E}_\psi\left(\overline{A}\sigma;\overline{W};\overline{I};\overline{U};\Phi\right) = 
	\left\|\overline{A}\overline{\varphi}\right\|_{\psi}.
\end{equation}
\end{theorem}

The next proposition deals with the optimal recovery problem of the sum of integral operators in the space $L_1(N,\nu)$.

\begin{cor}
\label{Cor2}
Let assumptions of Theorem~\ref{T2} hold. In addition, we assume that kernels $K_{ij}$ are $\nu_i\times\mu_j$-integrable ($i=1,\dots, l$, and $j=1,\dots, m$), $Y_i = L_1\left(N_i,\nu_i\right)$, and $\psi$ is $\ell_1$-norm on $\mathbb{R}^l$. Then
\[
	\begin{array}{l}
		\displaystyle\mathcal{E}_\psi\left(\overline{A}\sigma;\overline{W};\overline{I};\overline{U}\right) = \displaystyle \mathcal{E}_\psi\left(\overline{A}\sigma;\overline{W};\overline{I};\overline{U};\Phi\right)\\
		\qquad\qquad\qquad = \displaystyle\sum\limits_{i=1}^l\sum_{j=1}^{m} \int_{N_i}\int_{M_j'} K_{ij}\left(y,x\right) \varphi_j\left(x\right) \,d\mu_j(x)\,d\nu_i(y).
	\end{array}
\]
\end{cor}

\section{Applications} \label{S4}

In this section we demonstrate how main results of this paper can be applied to the problems of optimal recovery of the solutions to integral equations, and boundary and initial value problems for differential equations. 

In Section 4.1 and 4.2 we present optimal methods and error of recovery of solutions for Volterra and Fredholm integral equations of the second kind. 

We then present optimal methods and errors of recovery for solutions of systems of linear first order differential equations with constant coefficients, Poisson's equation, the heat and wave equations. Certainly, our approach is not restricted to optimal recovery of the solutions to mentioned equations and is applicable to a wider range of integral equations, ODE's, and PDE's.

\subsection{Linear Volterra integral equations of the second kind}

Let $[a,b]\subset \mathbb{R}$ be the interval, $\mu$ be the Lebesgue measure on $[a,b]$, function $f:[a,b]\to\mathbb{R}$ and kernel $k:[a,b]^2 \to \mathbb{R}$ be given, $x:[a,b]\to\mathbb{R}$ be unknown function. A linear Volterra equation of the second kind is the equation 
\begin{equation}
\label{Volterra}
	x(t) = f(t) + \int_a^t k(t,s)\,x(s)\,d\mu(s),\qquad t\in[a,b].
\end{equation}
We let $\Gamma$ be {\it the resolvent kernel} for $k$:
\[
	\Gamma(t,s) := \sum\limits_{n=1}^{\infty} k_n(t,s),\qquad t,s\in[a,b],
\]
where $k_1 = k$, and, for $n=2,3,\ldots$,
\[
	k_n(t,s) := \int_s^t k(t,\tau) k_{n-1}(\tau,s)\,d\mu(\tau),\qquad t,s\in[a,b].
\]

It is well known (see, for instance, ~\cite[Theorem~3.3]{Linz}) that for continuous function $f$ and kernel $k$ the solution to~(\ref{Volterra}) exists, is unique, and can be written in the form 
\begin{equation}
\label{Volterra_solution}
	x(t) = f(t) + \int_a^t \Gamma(t,s)\,f(s)\,d\mu(s),\qquad t\in[a,b].
\end{equation}

Let $\omega$ be the modulus of continuity, $Q$ be the set of $n\in\mathbb{N}$ points on $[a,b]$, and $e\in\mathbb{R}^n_+$. Let us consider the problem of optimal recovery of the solution to equation~(\ref{Volterra}) under assumptions that the values of function $f\in H^\omega_\mu([a,b])$ are known at the system of points $Q$ with $U_e$-error. 

By~(\ref{Volterra_solution}), the solution $x$ to problem~(\ref{Volterra}) can be considered as the image $x = Af$ of the function $f$ under the operator $A$, which is the sum of identity operator and an ntegral operator with the kernel $K(t,s) = k(t,s)\cdot \chi_{(a,t)}(s)$.

Let $Y\subset \mathfrak{M}([a,b],\mu)$ be a normed lattice containing the set of $B_\mu([a,b])$ of $\mu$-essentially bounded on $[a,b]$ functions. Then in view of Proposition~\ref{Upper_Estimate_for_Positive_Operator}, there holds true the following

\begin{cor}
	\label{Volterra_corollary}
	If the kernel $k$ is non-negative on $[a,b]^2$ then the method $\Phi = A\,L_{\omega,Q,e}$ is the optimal method of recovery of operator $A$ on the class $H^\omega_\mu([a,b])$ based on information $I_Q$ with $U_e$-error, and, furthermore,
	\[
	\begin{array}{rcl}
	\mathcal{E}\left(A;H^\omega_\mu([a,b]);I_Q;U_e\right) & = & \mathcal{E}\left(A;H^\omega_\mu([a,b]);I_Q;U_e;\Phi\right) = \left\|A \tau_{\omega,Q,e}\right\|_Y \\
	& = & \displaystyle \left\|\tau_{\omega,Q,e}(\cdot) + \int_{a}^{(\cdot)}\Gamma((\cdot),s)\,\tau_{\omega,Q,e}(s)\,d\mu(s)\right\|_Y.
	\end{array}
	\]
\end{cor}

It follows directly from Corollary~\ref{Volterra_corollary} that for $Y = L_1([a,b],\mu)$, 
\[
	\mathcal{E}\left(A;H^\omega_\mu([a,b]);I_Q;U_e\right) = \int_a^b \left(1 + \int_s^b \Gamma(t,s)\,d\mu(t)\right)\tau_{\omega,Q,e}(s)\,d\mu(s).
\]

\subsection{Linear Fredholm integral equations of the second kind}

Let $[a,b]\subset \mathbb{R}$, $\mu$ be the Lebesgue measure on $[a,b]$, $f:[a,b]\to\mathbb{R}$ be continuous function, $x:[a,b]\to\mathbb{R}$ be unknown function, and kernel $k:[a,b]^2\to\mathbb{R}$ be such that 
\begin{equation}
\label{square_integrability}
	\int_a^b \int_a^b |k(t,s)|^2\,d\mu(t)\,d\mu(s) < 1.
\end{equation}
The linear Fredholm integral equation of the second kind is the equation 
\begin{equation}
\label{Fredholm}
	x(t) = f(t) + \int_a^b k(t,s)\,x(s)\,d\mu(s).
\end{equation}
By $\Gamma$ we denote {\it the resolvent kernel} for $k$:
\[
	\Gamma(t,s) := \sum\limits_{n=1}^\infty k_n(t,s),\qquad t,s\in[a,b],
\]
where $k_1 = k$, and, for $n=2,3,\ldots$,
\[
	k_n(t,s) = \int_a^b k_{n-1}(t,u)\,k(u,s)\,d\mu(u),\qquad  t,s\in[a,b].
\]
It is well known (see, for instance, ~\cite[p.~44]{Corduneanu}) that the unique solution to~(\ref{Fredholm}) is given by 
\begin{equation}
\label{Fredholm_solution}
	x(t) = f(t) + \int_a^b \Gamma(t,s)\,f(s)\,d\mu(s),\qquad t\in[a,b].
\end{equation}

We let $\omega$ be a modulus of continuity, $Q$ be a set of $n\in\mathbb{N}$ points on $[a,b]$, $e\in\mathbb{R}^n_+$. We consider the problem of optimal recovery of the solution to equation~(\ref{Fredholm}) under assumption that the values of function $f\in H^\omega_\mu([a,b])$ are known at the system of points $Q$ with $U_e$-error. 

In virtue of~(\ref{Fredholm_solution}), the solution $x$ to problem~(\ref{Cauchy_problem_for_the_heat}) can be considered as the image $x = Af$ of the function $f$ under the operator $A$, which is the sum of identity operator and an integral operator with the kernel $K(t,s) = k(t,s)$.

Let $Y\subset\mathfrak{M}([a,b],\mu)$ be a normed lattice containing the space $B_\mu([a,b])$. Then by Proposition~\ref{Upper_Estimate_for_Positive_Operator}, there holds true the following

\begin{cor}
	\label{Fredholm_corollary}
	If the kernel $k$ is $\mu$-a.e. non-negative on $[a,b]^2$ and satisfies~(\ref{square_integrability}), then the method $\Phi = A\,L_{\omega,Q,e}$ is the optimal method of recovery of operator $A$ on the class $H^\omega_\mu([a,b])$ based on information $I_Q$ with $U_e$-error, and, furthermore,
	\[
	\begin{array}{rcl}
	\mathcal{E}\left(A;H^\omega_\mu([a,b]);I_Q;U_e\right) & = & \mathcal{E}\left(A;H^\omega_\mu([a,b]);I_Q;U_e;\Phi\right) = \left\|A \tau_{\omega,Q,e}\right\|_Y \\
	& = & \displaystyle \left\|\tau_{\omega,Q,e}(\cdot) + \int_{a}^{b}\Gamma((\cdot),s)\,\tau_{\omega,Q,e}(s)\,d\mu(s)\right\|_Y.
	\end{array}
	\]
\end{cor}

It follows directly from Corollary~\ref{Fredholm_corollary} that for $Y = L_1([a,b],\mu)$
\[
	\mathcal{E}\left(A;H^\omega_\mu([a,b]);I_Q;U_e\right) = \int_a^b \left(1 + \int_a^b \Gamma(t,s)\,d\mu(t)\right)\tau_{\omega,Q,e}(s)\,d\mu(s).
\]

\subsection{Optimal recovery of solutions to the systems of differential equations}

Next, we consider the problem of optimal recovery of the solution to initial value problem for the system of linear first order differential equations with constant coefficients. Let us introduce several notation: let $d\in\mathbb{N}$, $S$ be $d\times d$ matrix with real entries, $[a,b]\subset\mathbb{R}$ be a finite interval, $\mu$ be the Lebesgue measure on $[a,b]$. By $\overline{x}$ we denote vector function $\overline{x} = \left(x_1,\ldots,x_d\right)^T$ consisting of functions $x_i\in\mathfrak{M}\left([a,b],\mu\right)$, $i=1,\dots,d$. Finally, we let $\overline{q}$ be a continuous function, and $p\in\mathbb{R}^d$ be some point.

The system of linear nonhomogeneous equations has the form: 
\begin{equation}
\label{linear_differential_system}
	\left\{	
	\begin{array}{ll}
	\overline{x}'(t) = S\,\overline{x}(t) + \overline{q}(t), & t\in[a,b],\\
	\overline{x}\left(a\right) = p.& 
	\end{array}
	\right.
\end{equation}
It is well known that the solution to~(\ref{linear_differential_system}) exists, is unique, and is provided by:
\begin{equation}
\label{linear_solution}
	\overline{x}(t) = e^{S(t-a)}\,p + \int_a^{t} e^{S(u-a)}\,\overline{q}(u)\,d\mu(u),\qquad t\in[a,b],
\end{equation}
where $e^M$ stands for {\it the exponent} of matrix $M$ which is the series
\[
	e^{M} := \sum\limits_{j=0}^{\infty} \frac{M^j}{j!}.
\]

Next, we let $\omega_1,\ldots,\omega_d$ be given moduli of continuity, $W_1 := \mathbb{R}^d$ and 
\[
	W_2 := H^{\omega_1}_\mu([a,b])\times \ldots \times H^{\omega_d}_\mu([a,b]),
\]
$Q$ be the given set of $n\in\mathbb{N}$ points on $[a,b]$, $I_1 := id_{\mathbb{R}^d}:\mathbb{R}^d \to \mathbb{R}^d$ and $I_2 := \textrm{diag}\left(I_Q,\ldots,I_Q\right)$, consisting of $d$ operators $I_Q : C_\mu([a,b]) \to \mathbb{R}^n$, be information operators, $e\in\mathbb{R}^d_+$ and $e_1,\ldots,e_d\in \mathbb{R}^n_+$ be the errors describing information, $U_1 := U_{e}$ and $U_2 := U_{e_1}\times\ldots\times U_{e_d}$. 

Let us consider the problem of optimal recovery of the solution to the system~(\ref{linear_differential_system}) under assumptions that the initial value $p$ is known with $U_1$-error, and, for $i=1,\ldots,d$, the values of component $q_i$ of the function $\overline{q}\in W_2$ at the system of points $Q$ are known with $U_{e_i}$-error.

Next, for $i,j=1,\ldots,d$, by $k_{ij}(t)$, $t\in[a,b]$ we denote the element of matrix $e^{S(t-a)}$ located in the $i$th row and the $j$th column. In addition, we consider operators $B_{ij}:\mathbb{R} \to C_\mu([a,b])$ and $C_{ij}: B_\mu([a,b]) \to C_\mu([a,b])$ defined as follows
\[
	\left(B_{ij} x\right) (t) = k_{ij}(t) \cdot x,\qquad t\in[a,b],\quad x\in\mathbb{R},
\]
and
\[
	\left(C_{ij} x\right)(t) = \int_a^t k_{ij}(u)\,x(u)\,d\mu(u),\qquad t\in [a,b],\quad x \in B_\mu([a,b]).
\]
In view of~(\ref{linear_solution}), the solution $\overline{x}$ to~(\ref{linear_differential_system}) is the sum $\overline{x} = \overline{B} p + \overline{C} \overline{q}$ of images of initial value $p$ and function $\overline{q}$ under operator matrices $\overline{B}$ and $\overline{C}$, respectively. This observation allows us to consider the problem of optimal recovery of the solution to equation~(\ref{linear_differential_system}) as the problem of optimal recovery of the matrix operator $\overline{A} = \left(\overline{B}\;\overline{C}\right)$ on the class $\overline{W} = W_1\times W_2$ based on information $\overline{I} = \textrm{diag}\left(I_1,I_2\right)$ with $\overline{U} = U_1\times U_2$-error. 

Next, we remark that a square matrix $S$ is called {\it essentially non-negative} if every non-diagonal entry of this matrix is non-negative. It is well known that for such matrix $S$, operators $B_{ij}$ and $C_{ij}$ (the entries of matrix operator $\overline{A}$) are positive.  

For $i = 1,\ldots,d$, we let $Y_i\subset\mathfrak{M}\left([a,b],\mu\right)$ be a normed lattice, containing space $B_\mu([a,b])$, and let $\psi$ be a monotone norm in $\mathbb{R}^d$. Let also $L_1 := id_{\mathbb{R}^d}$, and $L_2 := \textrm{diag}\left(L_{\omega_1,Q,e_1},\ldots,L_{\omega_d,Q,e_d}\right)$. In addition, we set $\tau_1 := e$ and 
\[
	\tau_2 := \left(\tau_{\omega_1,Q,e_1},\ldots,\tau_{\omega_d,Q,e_d}\right)^T.
\]
Applying Theorem~\ref{Optimality}, we obtain the following 

\begin{cor}
\label{linear_differential_systems_solution_homogeneous}
Let $S$ be essentially non-negative $d\times d$ matrix. Then $\Phi = \overline{A}\,\overline{L}$, $\overline{L} = \textrm{diag}\left(L_1,L_2\right)$, is the optimal method of recovery of operator $\overline{A}$ on the class $\overline{W}$ based on information $\overline{I}$ with $\overline{U}$-error, and, furthermore,
\[
	\begin{array}{l}
		\displaystyle\mathcal{E}_\psi\left(\overline{A}; \overline{W}; \overline{I};\overline{U}\right) = \mathcal{E}_\psi\left(\overline{A}; \overline{W}; \overline{I}; \overline{U};\Phi\right) =  \left\|\overline{A}\,\overline{\tau}\right\|_{\psi} = \left\|\overline{B}e + \overline{C} \tau_2\right\|_\psi,\quad \overline{\tau} = \tau_1\times\tau_2.
	\end{array}
\]
\end{cor}

One can easily adjust the arguments of this section to solve the problem of optimal recovery of the solutions to the system of linear homogeneous equations and to the system of linear nonhomogeneous equations with homogeneous initial values.

\subsection{Optimal recovery of the solution to the Dirichlet problem for Poisson's equation}

In this section, we consider the problem of optimal recovery of the solution to the Dirichlet problem for Poisson's equation. We let $d\in\mathbb{N}$, $|\cdot|$ denote the standard norm in Euclidean space $\mathbb{R}^d$, $\mu$ be the standard Lebesgue measure in $\mathbb{R}^d$, $\Omega\subset \mathbb{R}^d$ be a bounded domain with $C^1$-boundary, $\sigma$ be the surface area measure on the boundary $\partial\Omega$. 
As usual, $\Delta$ stands for the Laplace operator. In addition, we let $C(\Omega\cup\partial\Omega)$ be the space of continuous on $\Omega\cup\partial\Omega$ real-valued functions, and $C^2(\Omega)$ be the space of continuous functions having continuous first, and second order partial derivatives inside $\Omega$.

The Dirichlet problem for Poisson's equation consists of finding a function $x\in C(\Omega\cup\partial\Omega)\cap C^2(\Omega)$, called the solution, which satisfies \begin{equation}
\label{Dirichlet_Poisson}
	\left\{\begin{array}{ll}
		-\Delta x(t) = f(t), & \overline{t}\in \Omega,\\
		x(t) = g(t), & \overline{t}\in \partial \Omega.
	\end{array}\right.
\end{equation}
It is well known (see~\cite{Evans}) that the solution $x$ to~(\ref{Dirichlet_Poisson}) exists, is unique, and can be presented in the form
\begin{equation}
\label{Dirichlet_Poisson_solution}
	x(t) = \int_{\Omega} G(t, s)\,f(s)\,d\mu(s) - \int_{\partial\Omega} \frac{\partial G}{\partial\bar{n}}(t,s)g(s)\,d\sigma(s),\quad t\in \Omega,
\end{equation}
where $G(t,s)$ is Green's function of the domain $\Omega$, and $\frac{\partial G}{\partial \bar{n}}(t,s)$ is the outer normal derivative of $G$.

Below, we assume that the domain $\Omega$ and it boundary $\partial \Omega$ are endowed with the respective Euclidean metric, and the metric which agrees with the surface area measure $\sigma$ on $\partial \Omega$. We consider classes $W_1 := H^{\omega_1}_{\mu} (\Omega)$ and $W_2 := H^{\omega_2}_{\sigma}(\partial\Omega)$, where $\omega_1,\omega_2$ are given moduli of continuity, finite sets of points $Q_1\subset \Omega$ and $Q_2\subset \partial \Omega$ consisting of, respectively, $n_1\in\mathbb{N}$ and $n_2\in\mathbb{N}$ points. We also assume that information operators $I_1 := I_{Q_1}$ and $I_2 := I_{Q_2}$, and sets $U_1 := U_{e_1}$ and $U_2 := U_{e_2}$, where $e_1\in\mathbb{R}^{n_1}_+$ and $e_2\in\mathbb{R}^{n_2}_+$, describing the error of information,  are given.

Let us consider the problem of optimal recovery of the solution to the problem~(\ref{Dirichlet_Poisson}) under assumptions that the values of functions $f\in W_1$ and $g\in W_2$ at systems of points $Q_1$ and $Q_2$ are known, respectively, with $U_{1}$ and $U_{2}$-errors.

By~(\ref{Dirichlet_Poisson_solution}), the solution $x$ to~(\ref{Dirichlet_Poisson}) is the sum $x = A_1f + A_2g$ of images of functions $f$ and $g$ under integral operators $A_1:B_\mu(\Omega) \to \mathfrak{M}\left(\Omega\cup\partial \Omega,\mu\right)$ and $A_2:B_\sigma(\partial \Omega) \to \mathfrak{M}\left(\Omega\cup\partial \Omega,\mu\right)$ respectively with kernels 
\[
	K_1(t,s) = G(t,s), \qquad t,s\in\Omega,
\]
and
\[
	K_2(t,s) = -\frac{\partial G}{\partial \bar{n}}(t,s),\qquad t\in \Omega,\quad s\in\partial\Omega.
\] 
Hence, the problem of optimal recovery of the solution to the problem~(\ref{Dirichlet_Poisson}) can be reformulated as the problem of optimal recovery of the matrix operator $\overline{A} = \left(A_1\;A_2\right)$ on the class $\overline{W} = W_1\times W_2$ based on information $\overline{I} = \textrm{diag}\left(I_1,I_2\right)$ with $\overline{U} = U_1\times U_2$-error. 

Since both operators $A_1$ and $A_2$ are positive, the assumptions of Theorem~\ref{T2} are satisfied. For convenience, we let $Y\subset \mathfrak{M}\left(\Omega\cup\partial\Omega,\mu\right)$ be a normed lattice, containing space $B_\mu(\Omega\cup\partial\Omega)$, $L_i := L_{\omega_i,Q_i,e_i}$, and $\tau_i := \tau_{\omega_i,Q_i,e_i}$, $i=1,2$. 

\begin{cor}
	\label{Dirichlet_Poisson_result}
	The operator $\Phi = \overline{A}\,\overline{L}$, with $\overline{L} = \textrm{diag}\left(L_1,L_2\right)$, is the optimal method of recovery of operator $\overline{A}$ on the class $\overline{W}$ based on information $\overline{I}$ with $\overline{U}$-error. Moreover, the optimal error is
	\[
		\mathcal{E}\left(\overline{A}; \overline{W}; \overline{I};\overline{U}\right) = \mathcal{E}\left(\overline{A}; \overline{W}; \overline{I};\overline{U};\Phi\right) = \left\|\overline{A}\,\overline{\tau}\right\|_Y, \qquad \overline{\tau} = \tau_1\times \tau_2.
	\]
\end{cor}

It follows directly from Corollary~\ref{Dirichlet_Poisson_result}, that for $Y = L_1\left(\Omega\cup\partial\Omega,\mu\right)$
\[
	\begin{array}{rcl}
		\mathcal{E}\left(\overline{A};\overline{W};\overline{I};\overline{U}\right) & = & \displaystyle \int_\Omega \tau_{1}(s)\int_{\Omega} G(t,s)\,d\mu(t)\,d\mu(s) \\ 
		& & \quad \displaystyle- \int_{\partial\Omega} \tau_{2}(s)\int_{\Omega} \frac{\partial G}{\partial \bar{n}} (t,s)\,d\mu(s)\,d\sigma(s).
	\end{array}
\]
In particular, when $\Omega$ is the disk of radius $r$ centered at the point $a\in\mathbb{R}^d$, we have
\[
	\mathcal{E}\left(\overline{A};\overline{W};\overline{I};\overline{U}\right) = \frac{1}{2}\int_{\Omega} \left(r^2 - |s-a|^2\right)\,\tau_{1}(s)\,d\mu(s) + \frac{r}{2} \int_{\partial \Omega} \tau_{2}(s)\,d\sigma(s). 
\]

One can apply similar arguments to solve the problem of optimal recovery of solutions to the Dirichlet problem for Laplace's equation, and to the homogeneous Dirichlet problem for Poisson's equation.

\begin{figure}[h!]
	\centering
	\includegraphics[angle=0, width=3.5in]{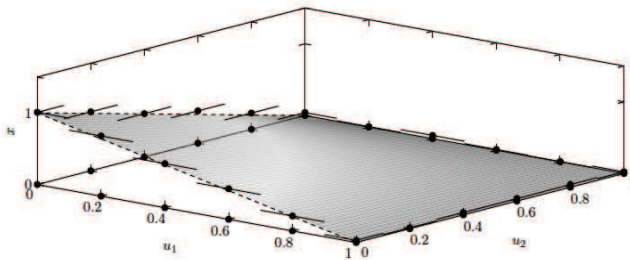}\label{pic4}
	\caption{Solution to the Laplace equation}
	\end{figure}

\begin{figure}[h!]
	\centering
	\includegraphics[angle=0, width=3.5in]{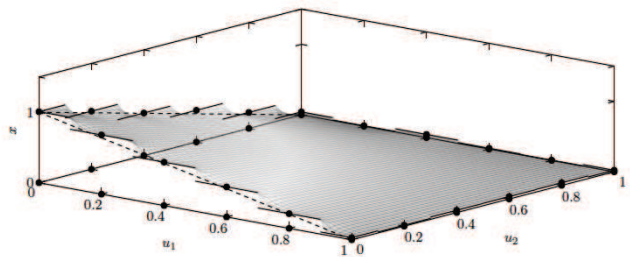}\label{pic5}
		\caption{Recovery of the solution to the Laplace equation}
	\end{figure}


\subsection{Optimal recovery of the solution to the initial value problem for heat equation}

The heat (or diffusion) equation describes the evolution in time of the density of some quantity such as heat, chemical concentration, etc. In the present section we apply results on optimal recovery of integral operators to the problem of optimal recovery of the solution to  initial value problems for this type of PDE.

In what follows, we use the following notation. Let $D = \mathbb{R}^d\times (0,+\infty)$, $d\in\mathbb{N}$. Let also $|\cdot|$ stand for the Euclidean norm in $\mathbb{R}^d$. For a function $x:D\to\mathbb{R}$, we denote by $\Delta x$ its space-coordinates Laplace operator, {\it i.e.}
\[
	\Delta x(u,t) := \sum\limits_{j=1}^{d} \frac{\partial^2 x}{\partial u_j^2}(u,t).
\]
In addition, we let $N$ to be a manifold in $D$. By $\mu$ we denote the standard Lebesgue measure corresponding to the dimension of the manifold, equipped with this measure, {\it i.e.} for the space $\mathfrak{M}\left(D,\mu\right)$ measure $\mu$ is $(d+1)$-dimensional Lebesgue measure, while for the space $\mathfrak{M}\left(\mathbb{R}^d,\mu\right)$ the measure $\mu$ is $d$-dimensional Lebesgue measure. We take $f:D\to\mathbb{R}$ to be a continuous, compactly supported function with continuous partial derivatives $\frac{\partial f}{\partial t}$, $\frac{\partial f}{\partial u_j}$, and $\frac{\partial^2 f}{\partial u_j \partial u_k}$, $j,k=1,\ldots,d$. Finally, let $g:\mathbb{R}^d\to\mathbb{R}$ be a continuous, bounded function. 

The Cauchy problem for the heat equation consists of finding a function $x\in C(D\cup\partial D)\cap C^2(D)$, called the solution, which satisfies 
\begin{equation}
\label{Cauchy_problem_for_the_heat}
	\left\{\begin{array}{ll}
		\frac{\partial x}{\partial t}(u,t) - \Delta x(u,t) = f(u,t), & (u,t)\in D,\\
		x(v,0) = g(v), & v\in\mathbb{R}^d.
	\end{array}\right.
\end{equation}
It is well known (see, for instance,~\cite{Ferretti}) that the solution $x$ to~(\ref{Cauchy_problem_for_the_heat}), satisfying a growth condition
\begin{equation}
\label{growth_condition}
	|x(u,t)| \leqslant \alpha e^{\beta |u|^2},\qquad (u,t)\in D,
\end{equation}
for some constants $\alpha,\beta>0$, exists, is unique, and can be presented in the form
\begin{equation}
\label{Cauchy_problem_for_the_heat_solution}
	\begin{array}{rcl}
	x(u,t) & = & \displaystyle \int_0^t\int_{\mathbb{R}^d} \frac{\exp{\left(-\frac{|u - v|^2}{4(t-s)}\right)}}{(4\pi(t-s))^{d/2}}\,f(v,s)\,d\mu(v,s) \\
	& & \displaystyle \quad + \int_{\RR^d} \frac{\exp{\left(-\frac{|u-v|^2}{4t}\right)}}{(4\pi t)^{d/2}}\,g(v)\,d\mu(v),\quad (u,t)\in D.
	\end{array}
\end{equation}

Below, we assume that $D$ and $\mathbb{R}^d$ are endowed with the Euclidean distance, $M_1\subset D$ and $M_2\subset\mathbb{R}^d$ are compact sets. We consider the classes $W_1 := \tilde{H}^{\omega_1}_{\mu}(M_1)$ and $W_2 := \tilde{H}^{\omega_2}_{\mu}(M_2)$, where $\omega_1,\omega_2$ are given moduli of continuity, finite sets of points $Q_1\subset M_1$ and $Q_2\subset M_2$ consisting of, respectively, $n_1\in\mathbb{N}$ and $n_2\in\mathbb{N}$ points. We also assume that information operators $I_1 := I_{Q_1}$ and $I_2 := I_{Q_2}$, and sets $U_1 := U_{e_1}$ and $U_2 := U_{e_2}$ describing the error of information, where $e_1\in\mathbb{R}^{n_1}_+$ and $e_2\in\mathbb{R}^{n_2}_+$, are given.

Let us consider the problem of optimal recovery of the solution to the problem~(\ref{Cauchy_problem_for_the_heat}) under assumptions that the values of functions $f\in W_1$ and $g\in W_2$ at systems of points $Q_1$ and $Q_2$ are known respectively with $U_1$ and $U_2$-errors.

Due to~(\ref{Cauchy_problem_for_the_heat_solution}), we can easily see that the solution $x$ to~(\ref{Cauchy_problem_for_the_heat}), satisfying the growth condition~(\ref{growth_condition}), is the sum $x = A_1f + A_2g$ of images of functions $f$ and $g$ under integral operators $A_1 : B_\mu(D) \to \mathfrak{M}(D,\mu)$ and $A_2:B_\mu(\mathbb{R}^d)\to \mathfrak{M}(D,\mu)$ with respective kernels
\[
	K_1((u,t),(v,s)) := \frac{\exp{\left(-\frac{|u - v|^2}{4(t-s)}\right)}}{(4\pi(t-s))^{d/2}}\cdot\chi_{(0,t)}(s),\qquad (u,t),(v,s)\in D,
\]
and
\[
	K_2((u,t),v) := \frac{\exp{\left(-\frac{|u - v|^2}{4t}\right)}}{(4\pi t)^{d/2}},\qquad (u,t)\in D,\quad v\in\mathbb{R}^d. 
\]
Hence, the problem of optimal recovery of the solution to problem~(\ref{Cauchy_problem_for_the_heat}) satisfying growth condition~(\ref{growth_condition}) can be reformulated as the problem of optimal recovery of the matrix operator $\overline{A} = \left(A_1\;A_2\right)$ on the class $\overline{W} = W_1\times W_2$ based on information $\overline{I} = \textrm{diag}\left(I_1,I_2\right)$ with $\overline{U} = U_1\times U_2$-error.

Since both operators $A_1$ and $A_2$ are positive, the assumptions of Theorem~\ref{T2} are satisfied. For convenience, we let $Y \subset \mathfrak{M}(N,\mu)$ be a normed lattice containing the space of $\mu$-essentially bounded and integrable on $N$ functions, $L_i = \tilde{L}_{\omega_i,Q_i,e_i}$ and $\tau_i = \tilde{\tau}_{\omega_i,Q_i,e_i}$, $i=1,2$.

\begin{cor}
	\label{Cauchy_the_heat}
	The operator $\Phi = \overline{A}\,\overline{L}$, $\overline{L} = \textrm{diag}\left(L_1,L_2\right)$, is the optimal method of recovery of operator $\overline{A}$ on the class $\overline{W}$ based on information $\overline{I}$ with $\overline{U}$-error. Moreover,
	\[
		\mathcal{E}\left(\overline{A};\overline{W};\overline{I};\overline{U}\right) = \mathcal{E}\left(\overline{A};\overline{W};\overline{I};\overline{U};\Phi\right) = \left\|\overline{A}\,\overline{\tau}\right\|_Y,\qquad \overline{\tau} = \tau_1\times\tau_2.
	\]
\end{cor}
In particular, for $Y = L_1(N,\mu)$ we obtain
\begin{enumerate}
	\item If $N = \mathbb{R}^d\times\left\{t_0\right\}$ where $t_0>0$, then 
	\[
	\displaystyle\mathcal{E}\left(\overline{A};\overline{W};\overline{I};\overline{U}\right) = \displaystyle \int_{\mathbb{R}^d}\int_0^{t_0} \tau_{1}\left(v,s\right)\,d\mu(s)\,d\mu(v) + \displaystyle \int_{\mathbb{R}^d} \tau_{2}\left(v\right)\,d\mu(v);
	\]
	\item If $N = \left\{u_0\right\}\times (0,+\infty)$ where $u_0\in\mathbb{R}^d$ is fixed, then
	\[
	\begin{array}{rcl}
	\displaystyle\mathcal{E}\left(\overline{A};\overline{W};\overline{I};\overline{U}\right) & = & \displaystyle \frac{4^{d-1}\Gamma\left(\frac d2 - 1\right)}{\pi^{d/2}}\left\{\int_0^{+\infty}\!\!\!\int_{\RR^d} \frac{\tau_{1}(v,s)}{\left|u_0 - v\right|^{d - 2}}\,d\mu(v)\,d\mu(s) \right. \\
	& & \quad + \displaystyle \left.\int_{\mathbb{R}^d} \frac{\tau_{2}\left(v\right)}{\left|u_0 - v\right|^{d - 2}}\,d\mu(v)\right\};
	\end{array}
	\]
	\item If $N = \left\{u_0\right\}\times \left\{t_0\right\}$ where $u_0\in\mathbb{R}^d$ and $t_0>0$ are fixed, then
	\[
	\begin{array}{l}
	\displaystyle\mathcal{E}\left(\overline{A};\overline{W};\overline{I};\overline{U}\right) = \displaystyle \displaystyle \int_0^{t_0}\int_{\RR^d} \frac{\exp{\left(-\frac{\left|u_0-v\right|^2}{4\left(t_0-s\right)}\right)} \tau_{1}(v,s)}{(4\pi\left(t_0-s\right))^{d/2}}\,d\mu(v)\,d\mu(s) \\
	\displaystyle\qquad\qquad\qquad\qquad\qquad+\int_{\mathbb{R}^d} \frac{\exp{\left(-\frac{\left|u_0-v\right|^2}{4t_0}\right)}}{\left(4\pi t_0\right)^{d/2}}\, \tau_{2}\left(v\right)\,d\mu(v).
	\end{array}
	\]
\end{enumerate}

Finally, we note that one can adjust the above arguments to solve the problem of optimal recovery of the solutions to the Cauchy problem for the homogeneous heat equation, and to the homogeneous Cauchy problem for the heat equation satisfying growth condition~(\ref{growth_condition}).

\subsection{Optimal recovery of the solution to initial value problem for the wave equation}

In this section we consider the problem of optimal recovery of the solution to the wave equation. We recall that the wave equation describes wave propagation in a media and is a simplified model for a vibrating string ($d=1$), membrane ($d=2$), or elastic solid ($d=3$). 

For the purpose of this section, we take $d=1,2,3$. Let $|\cdot|$ be the Euclidean norm in $\mathbb{R}^d$. We introduce notation of a $d$-dimensional ball $B_t(u)\subset \RR^d$ and a $d$-dimensional sphere $S_t(u)\subset \RR^{d-1}$ centered at point $u\in\mathbb{R}^d$ with radius $t>0$. In addition, we let $D = \mathbb{R}^d\times (0,+\infty)$. Similarly to the previous section, for a function $x:D\to\mathbb{R}$, we denote by $\Delta x$ its space-coordinates Laplace operator, let $N\subset D$ be a manifold, and $\mu$ be the standard Lebesgue measure corresponding to the dimension of the manifold, equipped with the measure. Let also $f:D\to\mathbb{R}$, $g:\mathbb{R}^d\to\mathbb{R}$, and $h:\mathbb{R}^d\to\mathbb{R}$ be some functions.

First, we let $d=1$ and consider the Cauchy problem for one-dimensional wave equation, which consists of finding a twice continuously differentiable function $x:D\to\mathbb{R}$, called a solution, satisfying the system of equations:
\begin{equation}
\label{wave_equation_1}
	\left\{\begin{array}{ll}
		x_{tt}(u,t) - x_{uu}(u,t) = f(u,t), & (u,t)\in D,\\
		x(u,0) = g(u), & u\in\mathbb{R},\\
		x_t(u,0) = h(u), & u\in\mathbb{R}.
	\end{array}\right.
\end{equation}
If $f$ is continuous, $g$ is twice continuously differentiable, and $h$ is continuously differentiable, then the unique solution to the problem~(\ref{wave_equation_1}) is delivered by the Dalambert formula:
\begin{equation}
\label{wave_equation_ihh_solution_1}
\begin{array}{rcl}
x(u,t) & = & \displaystyle \frac{1}{2}\int_0^t\!\! \int_{u-s}^{u+s} f(v,t-s)\,d\mu(v)\,d\mu(s) \\ 
& & \quad \displaystyle + \frac{g(u-t) + g(u+t)}{2} +  \displaystyle \frac{1}{2}\int_{u-t}^{u+t} h(v)\,d\mu(v).
\end{array}
\end{equation}
As we are interested in applications of Theorem~\ref{T2}, we would need to assume that $g$ and $h$ are compactly supported and have a majorant for their modulus of continuity. This means that $g$ and $h$ are continuous, but might be non-differentiable. Hence, the formula~(\ref{wave_equation_ihh_solution_1}) does not deliver the solution to the problem~(\ref{wave_equation_1}). Therefore, we follow~\cite{Petrovsky} to introduce the concept of the generalized solution to the problem~(\ref{wave_equation_1}). 

Let $g,h\in \tilde{C}_\mu(\Omega)$, where $\Omega$ is a compact in $\mathbb{R}$. Let also  $\left\{g^{(n)}\right\}_{n=1}^{\infty}$ be a sequence of twice continuously differentiable functions converging uniformly on $\Omega$ to the function $g$, as $n\to\infty$. Similarly, let $\left\{h^{(n)}\right\}_{n=1}^{\infty}$ be a sequence of continuously differentiable functions converging uniformly on $\Omega$ to $h$, as $n\to\infty$. The function $x(u,t)$ is called {\it the generalized solution} to the problem~(\ref{wave_equation_1}) if it is a limit of uniformly converging sequence of solutions $\left\{x^{(n)}(u,t)\right\}_{n=1}^{\infty}$ to the wave equation $x^{(n)}_{tt}(u,t) - x^{(n)}_{uu}(u,t) = f(u,t)$, $(u,t)\in D$, with initial conditions 
\[
	\begin{array}{c}
	x^{(n)}(u,0) = g^{(n)}(u),\qquad u\in\mathbb{R},\\
	x_t^{(n)}(u,0) = h^{(n)}(u),\qquad u\in\mathbb{R}.
	\end{array}
\]
Using the above definition, the Dalambert formula~(\ref{wave_equation_ihh_solution_1}) delivers the unique generalized solution to the problem~(\ref{wave_equation_1}).

Next we let $M_1\subset D$, and $M_2,M_3\subset \mathbb{R}$ be compact sets. For $i=1,2,3$, we let $\omega_i$ be modulus of continuity, $W_i := H^{\omega_i}_\mu(M_i)$ be a class of functions, 
$Q_i$ be the set of $n_i\in\mathbb{N}$ points on $M_i$, $I_i := I_{Q_i}$ be information operator, and $U_{e_i}$ be the error of information, $e_i\in\mathbb{R}_+^{n_i}$.

We see that the generalized solution $x$ to~(\ref{wave_equation_1}) can be presented as the sum $x = A_1f + A_2g + A_3h$ of images of functions $f$, $g$, and $h$ under linear positive operators $A_1: B_\mu(M_1) \to \mathfrak{M}(N,\mu)$, $A_2: B_\mu(M_2)\to \mathfrak{M}(N,\mu)$, and $A_3: B_\mu(M_3)\to \mathfrak{M}(N,\mu)$. This allows us to consider the problem of optimal recovery of the generalized solution $x$ to the wave equation~(\ref{wave_equation_1}) as the problem of optimal recovery of operator $\overline{A} = \left(A_1\;A_2\;A_3\right)$ on the class $\overline{W} = W_1\times W_2\times W_3$ based on information $\overline{I} = \textrm{diag}\,\left(I_1,I_2,I_3\right)$ with $\overline{U} = U_1\times U_2\times U_3$-error. 

Let us formulate the following corollary from Theorem~\ref{T1}. Let $Y\subset \mathfrak{M}(N,\mu)$ be the normed lattice containing essentially bounded and integrable on $N$ functions, and, for $i=1,2,3$, we let $L_i := L_{\omega_i,Q_i,e_i}$, $\tau_i := \tau_{\omega_i,Q_i,e_i}$. 

\begin{cor}
	\label{Cauchy_the_wave_1}
	The operator $\Phi = \overline{A}\,\overline{L}$, $\overline{L} = \textrm{diag}\left(L_1,L_2,L_3\right)$, is the optimal method of recovery of operator $\overline{A}$ on the class $\overline{W}$ based on information $\overline{I}$ with $\overline{U}$-error. Moreover,
	\[
	\mathcal{E}\left(\overline{A};\overline{W};\overline{I};\overline{U}\right) = \mathcal{E}\left(\overline{A};\overline{W};\overline{I};\overline{U};\Phi\right) = \left\|\overline{A}\,\overline{\tau}\right\|_Y,\qquad \overline{\tau} = \tau_1\times\tau_2\times \tau_3.
	\]
\end{cor}

In particular, for $Y = L_1(N,\nu)$, and $N = \mathbb{R}\times\{t_0\}$, where $t_0>0$ is fixed, we obtain
\[
	\begin{array}{rcl}
		\mathcal{E}\left(\overline{A};\overline{W};\overline{I};\overline{U}\right) & = & \displaystyle \int_0^{t_0} s \int_{-\infty}^{+\infty} \tau_1(v,t_0-s)\,d\mu(v)\,d\mu(s) + \int_{-\infty}^{\infty}\tau_2(v)\,d\mu(v) \\
		& & \qquad \displaystyle + t_0 \int_{-\infty}^{\infty} \tau_3(v)\,d\mu(v).
	\end{array}
\]

Next, we let $d=2,3$, and consider the Cauchy problem for the wave equation with zero initial form, {\it i.e.} $g\equiv 0$. This problem consists of finding a twice continuously differentiable function $x:D\to\mathbb{R}$, called a solution, which satisfies the system of equations:
\begin{equation}
\label{wave_equation_ihh}
\left\{\begin{array}{ll}
x_{tt}(u,t) - \Delta x(u,t) = f(u,t), & (u,t)\in D,\\
x(u,0) = 0, & u\in\mathbb{R}^d,\\
x_t(u,0) = h(u), & u\in\mathbb{R}^d.
\end{array}\right.
\end{equation}
If $f$ is continuous, and $h$ is twice continuously differentiable then the unique solution to the problem~(\ref{wave_equation_1}) is delivered by the Poisson formula when $d=2$:
\begin{equation}
\label{wave_equation_ihh_solution_2}
\begin{array}{rcl}
x(u,t) & = & \displaystyle \frac{1}{2\pi}\int_0^t\!\! \int_{B_{s}(u)} \frac{f(v,t-s)}{\sqrt{s^2 - |u-v|^2}}\,d\mu(v)\,d\mu(s)\\
& & \quad + \displaystyle \frac{1}{2\pi}\int_{B_{t}(u)}\frac{h(v)\,d\mu(v)}{\sqrt{t^2 - |u-v|^2}},
\end{array}
\end{equation}
and by the Kirchhoff formula for $d=3$:
\begin{equation}
\label{wave_equation_ihh_solution_3}
\begin{array}{rcl}
x(u,t) & = & \displaystyle \frac{1}{4\pi} \int_{B_{t}(u)} \frac{f\left(v,t - |u-v|\right)}{|u-v|}\,d\mu(v)\\ 
& & \quad + \displaystyle \frac{1}{4\pi t}\int_{S_{t}(u)} h(v)\,d\sigma(v),
\end{array}
\end{equation}
where $\sigma$ stands for the surface area measure of the sphere $S_{t}(u)$.

Similarly to the case $d=1$, we would need to assume that $h$ is compactly supported and has a majorant for the modulus of continuity. Therefore, we follow~\cite{Petrovsky}, and introduce the generalized solution to the problem~(\ref{wave_equation_ihh}) as follows. 

Let $h\in \tilde{C}_\mu(\Omega)$, where $\Omega$ is a compact in $\mathbb{R}$, and $\left\{h^{(n)}\right\}_{n=1}^{\infty}$ be the sequence of twice continuously differentiable functions converging uniformly on $\Omega$ to $h$, as $n\to\infty$. The function $x(u,t)$ is called {\it the generalized solution} to the problem~(\ref{wave_equation_ihh}) if it is a limit of uniformly converging sequence of solutions $\left\{x^{(n)}(u,t)\right\}_{n=1}^{\infty}$ to the wave equation $x^{(n)}_{tt}(u,t) - x^{(n)}_{uu}(u,t) = f(u,t)$, $(u,t)\in D$, with initial conditions 
\[
\begin{array}{c}
x^{(n)}(u,0) = 0,\qquad u\in\mathbb{R},\\
x_t^{(n)}(u,0) = h^{(n)}(u),\qquad u\in\mathbb{R}.
\end{array}
\]
Using the above definition, the Poisson and the Kirchhoff formulas deliver the unique generalized solution to problem~(\ref{wave_equation_ihh}).

Now, we let $M_1\subset D$ and $M_2 \subset \mathbb{R}^d$ be compact sets. For $i=1,2$, we let $\omega_i$ to be a modulus of continuity, $W_i := H^{\omega_i}_\mu(M_i)$ be a class of functions, $Q_i$ be the set of $n_i\in\mathbb{N}$ points on $M_i$, $I_i := I_{Q_i}$ be information operator, and $U_{e_i}$ be the error of information, $e_i\in\mathbb{R}_+^{n_i}$.

We see that the generalized solution $x$ to~(\ref{wave_equation_ihh}) can be presented as the sum $x = A_1f + A_2h$ of images of functions $f$ and $h$ under linear positive operators $A_1: B_\mu(M_1) \to \mathfrak{M}(N,\mu)$, and $A_2: B_\mu(M_2)\to \mathfrak{M}(N,\mu)$. This allows us to consider the problem of optimal recovery of the generalized solution $x$ to the wave equation~(\ref{wave_equation_ihh}) as the problem of optimal recovery of operator matrix $\overline{A} = \left(A_1\;A_2\right)$ on the class $\overline{W} = W_1\times W_2$ based on information $\overline{I} = \textrm{diag}\,\left(I_1,I_2\right)$ with $\overline{U} = U_1\times U_2$-error. 

From Theorem~\ref{T2} we obtain the following corollary. Let $Y\subset \mathfrak{M}(N,\mu)$ be the normed lattice containing $\mu$-essentially bounded and integrable on $N$ functions, and, for $i=1,2$, we let $L_i := L_{\omega_i,Q_i,e_i}$, $\tau_i := \tau_{\omega_i,Q_i,e_i}$.

\begin{cor}
	\label{Cauchy_the_wave_23}
	The operator $\Phi = \overline{A}\,\overline{L}$, $\overline{L} = \textrm{diag}\left(L_1,L_2\right)$, is the optimal method of recovery of operator $\overline{A}$ on the class $\overline{W}$ based on information $\overline{I}$ with $\overline{U}$-error. Moreover,
	\[
	\mathcal{E}\left(\overline{A};\overline{W};\overline{I};\overline{U}\right) = \mathcal{E}\left(\overline{A};\overline{W};\overline{I};\overline{U};\Phi\right) = \left\|\overline{A}\,\overline{\tau}\right\|_Y,\qquad \overline{\tau} = \tau_1\times \tau_2.
	\]
\end{cor}

In particular, for $Y = L_1(N,\mu)$, and $N = \mathbb{R}^d\times\left\{t_0\right\}$, where $t_0>0$ is fixed, we obtain
\[
	\mathcal{E}\left(\overline{A};\overline{W};\overline{I};\overline{U}\right) = \int_0^{t_0} s \int_{\mathbb{R}^d} \tau_1(v, t_0 - s)\,d\mu(v)\,d\mu(s) + t_0 \int_{\mathbb{R}^d} \tau_2(v)\,d\mu(v).
\]

\end{document}